# The harmonic $har$-index: an impromevement of the Hirsch $h$-index and the Egghe's $g$-index

Romeo Meštrović
Maritime Faculty Kotor, University of Montenegro,
Dobrota, 85330 Kotor, Montenegro
e-mail: romeo@ucg.ac.me

**Abstract**

In order to characterize the scientific output of scientists, in this paper we define the harmonic $har$-index whose values are positive integers. It is proved that $h \le har \le g$, where $h$ is the Hirsch index and $g$ is the Egghe's index. Despite the fact that the $har$-index is defined in a completely different way (including sum of reciprocals of citations of a researcher), based on our computatioanl results, we get the surprising fact that this index is highly correlated with the $hg$- index ($hg = \sqrt{hg}$) introduced by Alonso et al. (2010). Accordingly, we believe that the $har$-index keep the advantages of both measures as well as to minimize their disadvantages. In addition, it is much easier to calculate the values of $har$-index than those of $g$-index and $hg$-index.

## 1. Introduction

In the past two decades, numerous Scientometrics and Bibliometrics indicators were proposed to evaluate the scientific impact of individuals, institutions, colleges, universities and research teams. To appraise the scientific impact of scholars, institutions and research areas among others, several publication-based indicators are used, such as the size-dependent indicators (total number of citations and number of highly cited papers) and size-independent indicators (average number of citations per paper and proportion of highly cited publications) (Waltman, 2016), as well as citation frequency (life cycle) of papers etc. Based on the limitations of these indicators, Hirsch (2005) proposed a new indicator called *h*-index, whilst Egghe (2006a, 2006b, 2006c) defined and studied an improvement of the $h$-index, called the $g$-index.

The $h$-index is defined as the highest number of publications of a scientist that received $h$ or more citations (Hirsch, 2005). Hirsch (2007) reported results of an empirical study of the predictive power of the *h*-index compared with other bibliometric indicators. Hirsch indicated that the $h$-index is better than other indicators considered (total citation count, citations per paper and the total paper count) in predicting future scientific achievement. The dependence of $h$-index and $N_{cit}$ has been extensively studied by Egghe and Rousseau (2006), Glänzel (2006), Malesios (2015), Molinari and Molinari (2008) and Montazerian and al. (2019).

The $g$-index is defined as the highest rank such that the cumulative sum of the number of citations received is larger than or equal to the square of this rank. As noticed by Bihari et al. (2021), the $h$-index and the $g$-index got a lot of attention due to its simplicity, and several other indicators were proposed to extend the properties of the $h$-index and to overcome its shortcomings. Bornmann, Mutz, Hug and Daniel (2011) noticed that many studies have investigated many variants of the $h$-index, but these variants do not offer clear advantages over the original $h$-index (also see Anderson et al., 2008).

The $h$-index and its variants are extensively studied in Informetrics (both experimentally and theoretically). Numerous generalizations, extensions and variations of the $h$-index and the $g$-index were defined, studied and compared in the last two decades. In order to estimate the rang of the $h$-index as function of number of total citations of a researcher, Yong (2014) proposed the asymptotic expression for the $h$-index based on the asymptotic formula for the most likely Durfee square size for a partition in the set of all partitions of a positive integer $n$, established by Canfield, Corteel and Savage (1998) (also see Ekhad and Zeilberger, 2014 and Glänzel and Schubert, 1985). The converse question was quite recently considered by Meštrović and Dragović, B. (2023b).

The comprehensive recent review on $h$-index and its alternative indices was presented by Bihari et al. (2021) (see also Alonso et al., 2009).

In this paper, we define the harmonic $har$-index whose properties combine the advantages of the $h$-index and the $g$-index to produce a value near to the $h$-index than to $g$-index.

This article is organized as follows. Section 2 contains the definitions and notations of basic bibliometric indicators, as well as the definitions of $h$-, $g$-, $A$-, $R$-, $e$- and the $hg$-index. In Section 3 we define the new harmonic $har$-index and we prove that $h \leq har \leq g$. Some examples in which we compare related values of the $har$-index, $h$-index and the $g$-index are also presented in this section, as well as an analysis of the impact of increasing the number of researcher' citations on the increase his $har$-index. Using 14 Price awardees/experts in Bibliometrics data (bases on Scopus), the computational results and related discussion are presented in Section 4. Concluding remarks are given in Section 5.

Notice that in Appendix, there is a complete table of all 14 author's citations from Scopus database used in our calculations (Table A1), which is colored by the values of the bibliometric indices obtained in this paper. Finally, all our main computational results obtained in software Wolfram *Mathematica* 11 are presented.

## 2. The $h$-, $g$-, $A$-, $R$-, $e$ and $hg$-index

Suppose that a scientist $S$ has published $n$ publications whose number of citations $cit_1, cit_2, cit_3, ..., cit_n$ in considered database are usually ranked in decreased order, i.e.,

$$cit_1 \geq cit_2 \geq \cdots \geq cit_n \geq 0. \tag{1}$$

In this article, we will use the following notations and related notions.

$p$ - the total number of publications having at least one citation;

$N_{cit} = \sum_{j=1}^{p} cit_j$ - the total number of citations of a scientist in considered database;

$N_{cit}(s) = \sum_{j=1}^{p} cit_j$ - the total number of citations of a scientist in considered database up to the rank $s \in \{1, 2, ..., p\}$;

$N_{cit}(har) = \sum_{j=1}^{har} cit_j$ - the total number of citations of a scientist in considered database into $har$-core.

$N_{cit}(h) = \sum_{j=1}^{h} cit_j$ - the total number of citations of a scientist in considered database into $h$-core.

$N_{cit}(g) = \sum_{j=1}^{g} cit_j$ - the total number of citations of a scientist in considered database into $g$-core.

The Hirsch or $h$-*index*, introduced by Hirsch (2005) is defined as the highest rank $m$ such that the first $m$ publications of a scientist received $m$ or greater than $m$ citations in considered database, or equivalently,

$$h = \max\{m : cit_m \geq m, 1 \leq m \leq p\}. \quad (2)$$

The Egghe's $g$-*index* (Egghe, 2006a, b, c) is defined as the largest positive integer $k$ ($k = 1, 2, ..., n$) such that

$$\frac{1}{k}\sum_{j=1}^{k} cit_j \geq k. \quad (3)$$

From this definition we see that a scientist has $g$-index $g$ if $g$ is the highest rank such that his top $g$ publications have together, at least $g^2$ citations, while his top $g+1$ publications have, together less than $(g+1)^2$ citations. Zhang (2009a, 2009b) noticed that the $g$-index improves the $h$-index in that it not only keeps some advantages of the $h$-index but also counts citations from highly cited articles; however, the $g$-index suffers from a basic drawback. In fact, when $\sum_{j=1}^{n} cit_j > n^2$ the $g$-index has no definition.

Jin's $A$-*index* (the name was suggested by Rousseau, 2006), introduced by Jin (2006) is defined as the average number of citations received by the publications into $h$-core, i.e.,

$$A = \frac{1}{h}\sum_{j=1}^{h} cit_j, \quad (4)$$

which is under our notation equal to $N_{cit}(h)/h$.

Observe that the $A$-index achieves the same goal as the $g$-index, i.e. it takes into account the total number of citations included in the $h$-core. The $A$-index could potentially disadvantage prolific researchers, so the $A$-index was defined to remove this problem (Jin et al., 2007)

Another attempt to improve the insensitivity of the $h$-index to the number of citations to highly cited papers is the $R$-index, introduced in Jin et al. (2007. p. 857). The $R$-*index* as defined as

$$R = \sqrt{\sum_{j=1}^{h} cit_j}, \quad (5)$$

So as the $h$-index, also this measure takes into account the actual $cit_j$ - values in the $h$-core. Note that

$$R = \sqrt{hA} = \sqrt{N_{cit}(h)}. \qquad (6)$$

It is easy to show that (see Proposition 1 and Corollary in Jin et al, 2007)

$$h \leq g \leq A. \qquad (7)$$

Moreover, from $R = \sqrt{hA}$ and $A \geq h$ by Corollary in Jin et. al. (2007, p. 657) it follows that

$$h \leq R. \qquad (8)$$

Aimed at the same goal as the $g$-index, the $e$-index was proposed in 2009 by (Zhang, 2009a; also see Zhang, C.-T., 2013) with the aim of considering the contribution of excess citations, which are mostly from highly cited papers. It is defined as

$$e = \sqrt{\sum_{j=1}^{h} cit_j - h^2}. \qquad (9)$$

Clearly, one has $e = 0$ if and only if $cit_1 = cit_2 = \cdots = cit_h = h$. Since by the expression (5), $N_{cit}(h) = \sum_{j=1}^{h} cit_j = R^2 = hA$, we find that

$$e = \sqrt{hA - h^2} = \sqrt{R^2 - h^2} = \sqrt{N_{cit}(h) - h^2}. \qquad (10)$$

In 2010 Alonso et al. (2010) presented a new index, called the $hg$-index, is defined as the geometric mean of the $h$-index and the $g$-index, i.e.,

$$hg = \sqrt{hg}. \qquad (11)$$

Alonso et al. (2010) pointed out that this new index tries to fuse all the benefits of both previous measures and that tries to minimize the drawbacks that each one of them presented.

Obviously, one has $h \leq hg \leq g$, which together with the inequality (7) yields

$$h \leq hg \leq g \leq A. \qquad (12)$$

As already highlighted in some papers (see e.g., Jin et al., 2007 and Todeschini, 2011), the indices $g$-, $A$-, $R$- and $e$- are very sensitive to few highly-cited papers. In particular (see Todeschini, 2011), $g$-, $A$-, $R$- and $e$-indices appear to measure different characteristics of the author production because the information related to the number of citations prevails over that related to the number of papers in the $h$-core (e.g., the $g$-index). Burrell (2007) noticed that numerical investigations suggest that the $A$-index is a linear function of time and of the $h$-index, while the size of the Hirsch $h$-core has an approximate square-law relationship with time, and hence also with the $A$-index and the $h$-index. Notice that the $e$-index takes into account the contribution of excess citations.

### 3. The harmonic $har$-index

Under notions and related notations presented in Section 2, we give the following definition.

**Definition 1.** The *harmonic $har$-index* ie defined as the largest positive integer $k$ ($k = 1, 2, ..., p$) such that

$$\sum_{j=1}^{k} \frac{1}{cit_j} \leq 1, \quad (13)$$

i.e.,

$$har = \max \left\{ 1 \leq k \leq p : \sum_{j=1}^{k} \frac{1}{cit_j} \leq 1, \right\}, \quad (14)$$

where $p$ is a total number of publications of a researcher having at least one citation.

Notice that Definition 1 also holds if we replace $p$ by $n$ (which is a total number of publications of a researcher), assuming that $1/cit_j = 1/0 = \infty$, and hence, it holds

$$har \leq p. \quad (15)$$

**Theorem 1.** *Let $har$, $h$ and $g$ be the harmonic index, Hirsch index and Egghe's index, respectively. Then*

$$h \leq har \leq g. \quad (16)$$

*Proof.* First observe that $cit_i \geq h$ for all $i = 1, 2, ..., h$ and therefore,

$$\sum_{j=1}^{h} \frac{1}{cit_j} \leq \frac{h}{h} = 1. \quad (17)$$

The inequality (17) together with Definition 1.1 implies that

$$h \leq har. \quad (18)$$

By the well known Arithmetic mean-Harmonic mean inequality between $m$ ($m = 1, 2, ..$) positive real numbers $a_1, a_2, .., a_m$ (see, e.g., Mitrinović, 1970, p. 76), we have

$$H_m := \frac{m}{\frac{1}{a_1}+\frac{1}{a_2}+\cdots+\frac{1}{a_m}} \leq \frac{a_1+a_2+\cdots+a_m}{m}, \qquad (19)$$

where $H_m$ is the *Harmonic mean* and $A_m$ is the *Arithmetic mean* of numbers $a_1, a_2, ..., a_m$. The equality in (16) holds if and only if $a_1 = a_2 = \cdots = a_p$.

The inequality (19) immediately yields

$$a_1+a_2+\cdots+a_m \geq \frac{m^2}{\frac{1}{a_1}+\frac{1}{a_2}+\cdots+\frac{1}{a_m}}. \qquad (20)$$

Substituting $m = har$ into (20), and using the condition (13) from Definition 1 of *har*-index, we obtain

$$a_1+a_2+\cdots+a_{har} \geq (har)^2. \qquad (21)$$

The inequality (21) together with the definition of Egghe's $g$-index given by the inequality (3) yields

$$g \geq har. \qquad (22)$$

Finally, the inequlities (20) and (21) imply the inequality (15), which completes the proof of Theorem 1.

**Corollary 1.** *The following inequalities are true*:

$$h \leq har \leq g \leq A. \qquad (23)$$

*Proof.* The inequalities (23) immediately follows from the inequalities (7) and (16).

**Example 1.** To save the space, denote by $P = (cit_1, cit_2, ..., cit_p)$ the vector of all citations that are greater than equal to 1.

Consider the following four cases.

a)  $P = (5,4,3,2,1)$. Then $h = 3$, $har = 3$ (because of $1/5+1/4+1/3+1/2 = 1.283 > 1$, $1/5+1/4+1/3 = 0.7833 < 1$) and $g = 3$ (because of $5+4+3+2 = 14 < 16 = 4^2$ and $5+4+3 = 12 > 3^2$). Hence, $h = har = g = 3$.

b)  $P = (5,5,4,3,2,2)$. Then $h = 3$, $har = 4$ (because of $1/5+1/5+1/4+1/3+1/2 = 1.483 > 1$ $1/5+1/5+1/4+1/3 = 0.9833 < 1$) and $g = 4$ (because of $5+5+4+3+2 = 19 < 25 = 5^2$ and $5+5+4+2 = 16 = 6^2$). Hence, $3 = h < har = g = 4$.

c)  $P = (10,9,8,7,6,5,4,3,2,1)$. Then $h = 5$, $har = 6$ (because of $1/10+1/9+1/8+1/7+1/6+1/5+1/4 = 1.096 > 1$ and $1/10+1/9+1/8+1/7+1/6+1/5 = 0.846 < 1$) and $g = 7$ (because of $3+4+5+6+7+8+9+10 = 52 < 64 = 8^2$ and $4+5+6+7+8+9+10 = 49 = 49 = 7^2$). Hence, $5 = h < 6 = har < 7 = g$.

*d)* $P=(100,98,97,...,3,2,1)$. Then $h=50$, $har=63$ (because of $1/100+1/99+1/98+\cdots+1/39+1/39+1/38+1/37=1.012>1$ and $1/100+1/99+1/98+\cdots+1/39+1/39+1/38+1/37=1.012>1$ and $g=67$ (because of $34+35+36+\cdots+98+99+100=4489=67^2$ and $33+34+35+36+\cdots+98+99+100=4522<4624=68^2$). Hence, $har-h=13>4=g-har$, $har/h=1.26$ and $g/har=1.064$.

In Mathematics, for any positive integer $n$, the well known $n$-th harmonic number $H_n$ is defined as the sum of the reciprocals of the first $n$ positive integers, i.e.,

$$H_n = \sum_{j=1}^{n} \frac{1}{j}. \qquad (24)$$

This is the reason why we called the *har*-index from Definition 1 the harmonic index.

For given positive integer $m$, the *generalized $n$-th harmonic number of order $m$*, denoted as $H_n^{(m)}$, is defined as the sum of $m$-th powers of the reciprocals of the first $n$ positive integers, i.e.,

$$H_n^{(m)} = \sum_{j=1}^{n} \frac{1}{j^m}. \qquad (25)$$

Observe that $H_n^{(1)} = H_n$.

**Example 2.** For arbitrary fixed positive integers $n$ and $p$ such that $p \leq n$, define $cit_i = p-i+1$ for each $i=1,2,...,p$ and $cit_i = 0$ for each $i=p+1,...,n$. Then the $h$-index is defined as the largest integer $i$ such that $p-i+1 \geq i$. Hence, $h = \lfloor (p+1)/2 \rfloor$, where $\lfloor a \rfloor$ denotes the largest integer that does not exceed $a$.

By the estimate (3), the $g$-index is defined as the largest integer $k$ such that $\sum_{i=1}^{k}(p-i+1) \geq k^2$, whence it follows that $3k^2 - 2kp + k \leq 0$, i.e., $k(3k-2p-1) \leq 0$. Hence, $g = \lfloor (2p+1)/3 \rfloor$, and so $g/h = \lfloor (2p+1)/3 \rfloor / \lfloor (p+1)/2 \rfloor \approx 4/3 = 1.333$.

The harmonic *har*-index from this example can be approximately determined for large values of $p$. In this case, it is known (see, e.g., Miller et al., 2010) the following asymptotic expansion:

$$H_n := \sum_{j=1}^{n} \frac{1}{j} \approx \ln n + \gamma + \frac{1}{2n} - \frac{1}{12n^2} + o\left(\frac{1}{n^4}\right), \qquad (26)$$

where $\gamma \approx 0.5772156649$ is the Euler Mascheroni constant.

Then using the approximative expression (26) in the form $\sum_{j=1}^{n} \frac{1}{j} \approx \ln n + \gamma$, the inequality (13) from Definition 1 becomes

$$\sum_{j=1}^{k}\frac{1}{cit_j} = \sum_{j=1}^{k}\frac{1}{p-j+1} = \sum_{j=1}^{p}\frac{1}{j} - \sum_{i=1}^{p-k}\frac{1}{i} \approx (\ln p + \gamma) - (\ln(p-k) + \gamma) = \ln\frac{p}{p-k} \leq 1, \qquad (27)$$

whence we find that $k \leq \frac{e-1}{e}p = 0.6321205588285577p$. It follows that $har \approx \lfloor 0.6321205588285577p \rfloor$. Finally, for large values $p$ we find that

$$har/h \approx \lfloor 0.6321205588285577p \rfloor / \lfloor (p+1)/2 \rfloor \approx 1.2642411176571153$$

and

$$g/har \approx \lfloor 2(p+1)/3 \rfloor / \lfloor 0.6321205588285577p \rfloor \approx 1.0546511379128842.$$

Let's assume that the current harmonic index of a researchers is equal to *harm*. Now we will briefly analyze how the change of quote affects the sum $S$ defined as

$$S := \sum_{j=1}^{har}\frac{1}{cit_j}. \qquad (28)$$

Notice that by Definition 1, it holds

$$S \leq 1 \quad \text{and} \quad S' := S + \frac{1}{cit_{har+1}} > 1. \qquad (29)$$

We will called the set $H := \{cit_1, cit_2, \ldots, cit_{har}\}$ the *harm–core*.

If the citation $cit_j$ ($1 \leq j \leq har$) increases by 1, then the sum $S$ defined by (28) will decrease by

$$\delta_j := \frac{1}{cit_j} - \frac{1}{cit_j + 1} = \frac{1}{cit_j(cit_j+1)}. \qquad (30)$$

**Hence, if** $m$ citations $cit_{j_1}, cit_{j_2}, \ldots, cit_{j_m}$ ($1 \leq j_1, j_2, \ldots, j_m \leq har$ with the set $J := \{j_1, j_2, \ldots, j_m\}$) increase by $k_1, k_2, \ldots, k_m$ times, rexpectively, then the sum $S$ defined by (28) will decrease by

$$\delta := \sum_{j \in J} k_j \delta_j := \sum_{j \in J} \frac{k_j}{cit_j(cit_j+1)}. \qquad (31)$$

In this case, if it is satisfied the condition

$$har + \frac{1}{cit_{har+1}} - \sum_{j \in J} \frac{k_j}{cit_j(cit_j + 1)} \leq 1, \qquad (32)$$

then by Definition 1, the harmonic $har$-index will increase for at least 1.

As it is noticed in Abstract, based on our computatioanl results, we get the surprising fact that this index is highly correlated with the $hg$-index ($hg = \sqrt{hg}$) index introduced by Alonso et al. (2010).

From the above considerations we see that the increase of the number of citations of publications with fewer citations has a greater contribution to the increase of the value of $har$-index compared to the increase of the number of citations of publications with more citations. Increasing any citation in the $har$ core (which is wider than the $h$-core and narrower than the $g$-core) automatically decreases the sum of the reciprocal values of the citations, which increases the probability of an increase of the $har$-index. Acoordingly, the $har$-index can be considered as a dynamical type bibliometric indices. It should be pointed out that among all commonly used indices in bibliometric literature, the $h$-index, the $g$-index and the $har$-index are the only ones that have integer values.

As the $har$-index is highly correlated with the the $hg$-index, the following observation (the benefit for the the $hg$-index), written by Alonso et al. (2010, p. 395) applies to it to the the $har$-index: "The $hg$-index takes into account the cites of the highly cited papers (the $h$-index is insensitive to highly cited papers) but it significantly reduces the impact of single very high cited papers (a drawback of the $g$-index), thus achieving a better balance between the impact of the majority of the best papers of the author and very highly cited ones."

## 4. Computational results and discussion

Our computations are motivated by the fact that Glänzel and Persson (2005) calculated the $h$-index for the 14 Price awardees who were still active in quantitative studies of science (based on their published journal papers from January 1986 to August 2005 extracted from Web of Science (WoS) database). The researchers whose bibliometric indicators and the indices are calculated in this Section, using their Scopus profiles from February 21 2023 are the following 14 Price awardees/leading experts in Scientometrics:

Leot Leydesdorff (S), Wolfgang Glänzel (S), Henk Moed (S), Anthony Van Raan (S), Ronald Rousseau (S), Andras Schubert (S), Ben Martin (S), Francis Narin (S), Eugene Garfield (S), Tibor Braun (S), Henry Small (S), Leo Egghe (S), Peter Ingwersen (S) and Howard D. White (S).

Using the notions, notations and the expressions presented in Sections 2 and 3, data of Table A1 from Appendix (Table A2 in Meštrović and Dragović, 2023a) and Table 1 in in Meštrović and Dragović (2023a) we obtain the following Table 1.

**Table 1**
Price awardees data (bases on Scopus, February 21 2023). The bibliometric indicators and indices.

| | L. Leydesdorff | W. Glänzel | H.F. Moed | A.F.J Van Raan | R. Rousseau | A. Schubert | B. Martin |
|---|---|---|---|---|---|---|---|
| $p$ | 406 | 258 | 127 | 124 | 295 | 141 | 85 |
| $N_{cit}$ | 25005 | 11766 | 7606 | 8308 | 8053 | 7587 | 7597 |
| $h$ | 79 | 61 | 49 | 48 | 43 | 42 | 38 |
| $N_{cit}(h)$ | 17360 | 8049 | 6351 | 6833 | 5203 | 6402 | 7048 |
| $N_{cit}(h)/h$ | 219.747 | 131.95 | 129.612 | 142.235 | 121 | 151.143 | **185.48** |
| $e$ | 105.447 | 65.788 | 62.849 | 67.298 | 52.934 | 68.103 | 74.860 |
| $e/h$ | 1.335 | 1.078 | 1.283 | 1.402 | 1.231 | 1.622 | **1.970** |
| $g$ | 145 | 99 | 86 | 89 | 79 | 85 | 85 |
| $N_{cit}(g)$ | 21225 | 9810 | 7396 | 8060 | 6300 | 7350 | 7597 |
| $h/g$ | 0.546 | 0.616 | 0.570 | 0.539 | 0.544 | 0.494 | **0.447** |
| $hg := \sqrt{hg}$ | 107.028 | 77.711 | 64.915 | 65.361 | 58.284 | 59.750 | 56.833 |
| $A(h,g) := (h+g)/2$ | 112 | 80 | 67.5 | 68.5 | 61 | 63.5 | 61.5 |
| $H(h,g) := 2hg/(h+g)$ | 102.277 | 75.488 | 62.430 | 62.365 | 55.689 | 56.220 | 52.520 |
| $har$ | 106 | 80 | 64 | 66 | 57 | 56 | 50 |
| $hg/har$ | **1.0097** | **0.971** | **1.014** | **0.990** | **1.023** | **1.067** | **1.137** |
| $A(h,g)/har$ | 1.057 | **1.000** | 1.057 | 1.038 | 1.070 | 1.134 | 1.230 |
| $har/H(h,g)$ | **1.036** | **1.060** | **1.025** | **1.058** | **1.024** | **0.996** | **0.952** |
| $g/har$ | 1.368 | 1.238 | 1.344 | 1.349 | 1.386 | 1.518 | **1.700** |
| $har/h$ | 1.342 | 1.311 | 1.306 | 1.375 | 1.326 | 1.333 | 1.316 |
| $N_{cit}(har)$ | 19274 | 9052 | 6933 | 7544 | 5698 | 6848 | 7349 |
| $S_{har+1} := \sum_{j=1}^{har+1} \frac{1}{cit_j}$ | 1.009 | 1.001 | **1.00035** | 1.023 | 1.014 | 1.022 | 1.0268 |
| $S_{har} := \sum_{j=1}^{har} \frac{1}{cit_j}$ | 0.993 | 0.980 | 0.966 | 0.991 | 0.981 | 0.984 | 0.977 |
| $S_h := \sum_{j=1}^{h} \frac{1}{cit_j}$ | 0.610 | 0.618 | 0.572 | 0.510 | 0.636 | 0.536 | 0.486 |
| $S_g := \sum_{j=1}^{g} \frac{1}{cit_j}$ | 1.787 | 1.461 | 2.107 | 2.048 | 1.866 | 2.761 | 11.490 |

**Table 1 -Continued**
Price awardees data (bases on Scopus, February 2023). The bibliometric indicators and indices.

| | F. Narin | E. Garfield | T. Braun | H. Small | L. Egghe | P. Ingwersen | H.D. White |
|---|---|---|---|---|---|---|---|
| $p$ | 64 | 106 | 216 | 57 | 211 | 88 | 37 |
| $N_{cit}$ | 7209 | 11515 | 5680 | 7693 | 5640 | 3606 | 2399 |
| $h$ | 38 | 37 | 37 | 34 | 30 | 27 | 19 |
| $N_{cit}(h)$ | 6823 | 10509 | 3566 | 7471 | 3995 | 2952 | 2332 |
| $N_{cit}(h)/h$ | 179.553 | **284.027** | 153.514 | 219.735 | 133.167 | 109.333 | 122.737 |
| $e$ | 73.342 | 95.603 | 46.872 | 79.467 | 55.633 | 47.149 | 44.396 |
| $e/h$ | 1.930 | **2.584** | 1.267 | 2.337 | 1.854 | 1.746 | 2.337 |
| $g$ | 68 | 106 | 66 | 57 | 69 | 60 | 28 |
| $N_{cit}(g)$ | 7209 | 11359 | 4373 | 7693 | 4807 | 3516 | 2399 |
| $h/g$ | 0.559 | **0.350** | 0.561 | 0.596 | 0.435 | 0.450 | 0.679 |
| $hg := \sqrt{hg}$ | 50.833 | 62.626 | 49.417 | 44.023 | 45.497 | 40.249 | 23.065 |
| $A(h,g) := (h+g)/2$ | 53 | 71.5 | 51.5 | 45.5 | 49.5 | 43.5 | 23.5 |
| $H(h,g) := 2hg/(h+g)$ | 48.754 | 54.583 | 47.417 | 42.593 | 41.818 | 37.241 | 22.638 |
| $har$ | 49 | 51 | 50 | 43 | 42 | 38 | 23 |
| $hg/har$ | **1.037** | **1.228** | **0.988** | **1.024** | **1.083** | **1.059** | **1.0028** |
| $A(h,g)/har$ | 1.082 | 1.402 | 1.030 | 1.011 | 1.179 | 1.145 | 1.022 |
| $har/H(h,g)$ | **1.005** | **0.934** | **1.055** | **1.010** | **1.004** | **1.020** | **1.016** |
| $g/har$ | 1.388 | **2.078** | 1.320 | 1.326 | 1.643 | 1.579 | 1.217 |
| $har/h$ | 1.289 | 1.378 | 1.351 | 1.264 | 1.400 | 1.407 | 1.211 |
| $N_{cit}(har)$ | 7092 | 10838 | 3971 | 7633 | 4309 | 3220 | 2388 |
| $S_{har+1} := \sum_{j=1}^{har+1} \frac{1}{cit_j}$ | 1.015 | 1.026 | 1.037 | 1.027 | 1.025 | 1.002 | 1.027 |
| $S_{har} := \sum_{j=1}^{har} \frac{1}{cit_j}$ | 0.943 | 0.976 | **0.999693** | 0.928 | 0.981 | 0.956 | 0.827 |
| $S_h := \sum_{j=1}^{h} \frac{1}{cit_j}$ | 0.462 | 0.370 | 0.580 | 0.418 | 0.519 | 0.503 | 0.449 |
| $S_g := \sum_{j=1}^{g} \frac{1}{cit_j}$ | 2.950 | 7.574 | 1.595 | 7.916 | 2.471 | 2.737 | 4.027 |

From Table 1 we see that the estimate $|hg/har - 1| < 0.084$ is true for 12 scientists, except E. Garfieled ($hg/har = 1.228$) and B. Martin ($hg/har = 1.137$). Notice that E. Garfieled died in 2017 and he published his last paper in 2016. It seems that this is the reason why its values $N_{cit}(h)/h$ (284.027) and $e/h$ (2.584) are relatively large numbers. Because of this, its value $h/g$ (0.35) is small. Consequently, the E. Garfield's value $g/har$ is relatively very big and it is equal to 2.078. B. Martin's values $N_{cit}(h)/h$, $e/h$ and $h/g$ are equal 185.48, 1.970 and 0.447, respectively. Also B. Martin has the second highest value $g/har$ and it is equal to 1.700. Observe also that E. Garfield's and B. Martin's values of $1 - har/H(h,g)$ are 0.066 and 0.046, respectively.

Notice that the values $hg/har - 1$ of H.D. White and L. Leydesdorf are very close to zero and they are equal to 0.0028 and 0.0097, respectively. From this together with our additional computational results we can conclude that the $har$-index and the $hg$-index ($hg = \sqrt{hg}$) are highly correlated (it seems more correlated than the $g$-index and the $R$-index ($R = \sqrt{N_{cit}(h)}$). This is a surprising fact, taking into acount the fact that the harmonic $har$-index is defined in a completely different way than the $h$-index and the $g$-index.

Note also that the values $har/h$ from Table 1 vary between in 1.211 and 1.407. T. Braun has the value $S_{har}$ quite close to 1 and it is equal to 0.999693. The values $S_h$ vary between in 0.370 and 0.636, while he values $S_g$ vary between in 1.461 and 11.490.

Observe also that the estimate $har \leq hg$ is true for 11 scientists, except W. Glänzel, A.F.J Van Raan and T. Braun. However, the estimate $har \leq (h+g)/2$ is true for all 14 scientists, while for W. Glänzel, $har = (h+g)/2 = 80$. From Table 1 we also see that the estimate $|har/H(h,g) - 1| < 0.07$ is satisfied for all 14 scientists. From this together with our additional computational results we can conclude that the $har$-index and the $H(h,g)$-index are highly correlated ($H(h,g) := 2hg/(h+g)$ is the harmonic mean of the $h$-index and the $g$-index).

From Table 1 we also see that 14 considered scientists are ranked with respect to the $h$-index very close with respect to the $hg$-index and the $har$-index.

## 5. Concluding Remarks

In this article, we introduce a new bibliometric index, the so called harmonic $har$-index. Our computational results show that $har$-index is highly correlated with the $hg$-index ($hg = \sqrt{hg}$) introduced by Alonso et al. (2010). We believe that the $har$-index keep the advantages of both Hirsch $h$-index and Egghe's $g$-index, as well as to minimize their disadvantages, but it is natural that bibliometric experts will judge these two facts.

**Compilance with Ethical Standards**

• The author of this manuscript has no competing interests to declare that are relevant to the content of this manuscript.

# Appendix

**Table A1.** Price awardees data: the all citations in domain *D* of asymptotic normality of citations (bases on Scopus, 21st February 2023)

| | Authors | | | | | | | | | | | | | |
|---|---|---|---|---|---|---|---|---|---|---|---|---|---|---|
| | L. Leyde-sdorff | W. Glänzel | H. Moed | A. van Raan | R. Rousseau | A. Schubert | B. Martin | F. Narin | E. Garfield | T. Braun | H. Small | L. Egghe | P. Ingwe-Rsen | H.D. White |
| 1 | 3964 | 534 | 417 | 522 | 1050 | 1858 | 1828 | 803 | 1888 | 449 | 2846 | 1455 | 482 | 1029 |
| 2 | 647 | 449 | 417 | 455 | 484 | 449 | 728 | 780 | 1711 | 362 | 592 | 478 | 414 | 213 |
| 3 | 544 | 417 | 380 | 381 | 478 | 362 | 411 | 519 | 1635 | 255 | 540 | 242 | 348 | 178 |
| 4 | 479 | 364 | 350 | 350 | 265 | 273 | 391 | 481 | 597 | 158 | 410 | 206 | 232 | 148 |
| 5 | 424 | 273 | 303 | 303 | 242 | 255 | 371 | 443 | 589 | 151 | 242 | 153 | 166 | 124 |
| 6 | 422 | 255 | 290 | 290 | 175 | 211 | 362 | 354 | 430 | 147 | 197 | 124 | 156 | 123 |
| 7 | 410 | 249 | 288 | 284 | 153 | 159 | 282 | 339 | 264 | 138 | 193 | 110 | 155 | 104 |
| 8 | 401 | 220 | 248 | 283 | 126 | 158 | 244 | 306 | 255 | 135 | 190 | 107 | 133 | 72 |
| 9 | 376 | 211 | 197 | 276 | 118 | 155 | 223 | 283 | 252 | 127 | 186 | 88 | 113 | 54 |
| 10 | 352 | 201 | 183 | 206 | 111 | 151 | 196 | 226 | 238 | 120 | 157 | 81 | 70 | 50 |
| 11 | 352 | 188 | 182 | 172 | 110 | 147 | 179 | 193 | 238 | 119 | 154 | 71 | 57 | 35 |
| 12 | 301 | 163 | 174 | 166 | 97 | 133 | 131 | 183 | 178 | 96 | 151 | 70 | 56 | 32 |
| 13 | 291 | 159 | 166 | 157 | 96 | 127 | 128 | 159 | 155 | 77 | 144 | 70 | 49 | 29 |
| 14 | 266 | 155 | 152 | 152 | 88 | 120 | 122 | 151 | 154 | 75 | 129 | 68 | 46 | 28 |
| 15 | 252 | 147 | 146 | 150 | 87 | 119 | 114 | 138 | 153 | 73 | 113 | 65 | 45 | 25 |
| 16 | 242 | 145 | 122 | 146 | 75 | 105 | 101 | 121 | 151 | 69 | 112 | 50 | 42 | 24 |
| 17 | 241 | 141 | 118 | 142 | 74 | 98 | 77 | 114 | 139 | 67 | 111 | 48 | 42 | 23 |
| 18 | 232 | 135 | 118 | 136 | 72 | 92 | 73 | 94 | 132 | 61 | 102 | 47 | 41 | 21 |
| 19 | 225 | 133 | 117 | 130 | 70 | 89 | 72 | 87 | 114 | 59 | 98 | 46 | 39 | **20** |
| 20 | 221 | 127 | 103 | 111 | 67 | 75 | 69 | 81 | 113 | 57 | 90 | 45 | 39 | 19 |
| 21 | 184 | 126 | 95 | 110 | 66 | 73 | 68 | 76 | 111 | 55 | 81 | 42 | 38 | 17 |
| 22 | 179 | 126 | 92 | 103 | 64 | 73 | 66 | 73 | 91 | 55 | 79 | 41 | 35 | 15 |
| 23 | 171 | 121 | 88 | 98 | 63 | 72 | 65 | 73 | 88 | 52 | 71 | 40 | 35 | 5 |
| 24 | 166 | 111 | 80 | 94 | 61 | 71 | 64 | 68 | 83 | 52 | 57 | 40 | 33 | 5 |
| 25 | 165 | 109 | 76 | 92 | 58 | 70 | 63 | 67 | 74 | 51 | 51 | 39 | 29 | 2 |
| 26 | 157 | 108 | 75 | 86 | 56 | 69 | 58 | 61 | 73 | 47 | 49 | 35 | 29 | 2 |
| 27 | 154 | 96 | 75 | 82 | 56 | 69 | 57 | 61 | 69 | 46 | 47 | 34 | **28** | 1 |
| 28 | 147 | 94 | 73 | 81 | 55 | 64 | 55 | 55 | 65 | 45 | 44 | 34 | 27 | 1 |
| 29 | 147 | 92 | 71 | 81 | 51 | 63 | 50 | 54 | 64 | 45 | 43 | 34 | 27 | |
| 30 | 147 | 89 | 70 | 77 | 49 | 62 | 48 | 50 | 63 | 45 | 41 | **32** | 26 | |
| 31 | 141 | 89 | 69 | 76 | 48 | 62 | 48 | 43 | 63 | 43 | 40 | 30 | 25 | |
| 32 | 138 | 89 | 69 | 75 | 48 | 57 | 47 | 43 | 53 | 42 | 39 | 29 | 24 | |
| 33 | 135 | 85 | 63 | 73 | 47 | 53 | 46 | 43 | 51 | 40 | 38 | 29 | 24 | |
| 34 | 133 | 85 | 62 | 72 | 46 | 52 | 46 | 42 | 49 | 40 | 34 | 28 | 24 | |
| 35 | 133 | 84 | 61 | 72 | 45 | 47 | 45 | 41 | 45 | 39 | 21 | 27 | 23 | |
| 36 | 129 | 84 | 60 | 68 | 45 | 46 | 40 | 40 | 42 | 37 | 21 | 26 | 23 | |
| 37 | 129 | 83 | 60 | 65 | 45 | 46 | 40 | 40 | **39** | 37 | 20 | 25 | 23 | |
| 38 | 127 | 80 | 58 | 63 | 45 | 45 | **40** | 38 | 30 | 36 | 19 | 24 | 22 | |
| 39 | 127 | 80 | 57 | 62 | 44 | 45 | 32 | 37 | 28 | 34 | 18 | 24 | 22 | |
| 40 | 126 | 78 | 57 | 60 | 44 | 43 | 32 | 34 | 27 | 33 | 17 | 24 | 22 | |
| 41 | 125 | 76 | 56 | 59 | 43 | 42 | 31 | 32 | 27 | 33 | 17 | 24 | 18 | |
| 42 | 125 | 76 | 55 | 59 | 43 | **42** | 28 | 27 | 25 | 33 | 15 | 24 | 17 | |
| 43 | 125 | 75 | 53 | 59 | 43 | 40 | 24 | 24 | 23 | 32 | 14 | 23 | 16 | |
| 44 | 125 | 75 | 53 | 55 | 42 | 37 | 23 | 21 | 22 | 31 | 10 | 22 | 16 | |
| 45 | 124 | 75 | 52 | 52 | 42 | 36 | 23 | 21 | 22 | 31 | 9 | 22 | 15 | |
| 46 | 123 | 74 | 52 | 51 | 42 | 36 | 23 | 20 | 21 | 30 | 9 | 22 | 15 | |
| 47 | 122 | 73 | 50 | 48 | 41 | 36 | 22 | 18 | 21 | 29 | 6 | 22 | 14 | |
| 48 | 122 | 73 | 49 | **48** | 40 | 33 | 21 | 18 | **21** | 28 | 4 | 22 | 14 | |
| 49 | 118 | 73 | **49** | 47 | 39 | 32 | 21 | 17 | 21 | 28 | 4 | 20 | 13 | |
| 50 | 108 | 72 | 47 | 45 | 39 | 31 | 21 | 14 | 21 | 27 | 4 | 20 | 13 | |
| 51 | 106 | 72 | 45 | 41 | 39 | 31 | 20 | 12 | 20 | 27 | 4 | 20 | 13 | |
| 52 | 102 | 71 | 44 | 41 | 37 | 27 | 20 | 12 | 20 | 27 | 2 | 19 | 12 | |
| 53 | 100 | 71 | 43 | 39 | 37 | 27 | 20 | 11 | 20 | 27 | 1 | 19 | 12 | |
| 54 | 98 | 70 | 43 | 39 | 33 | 27 | 17 | 11 | 19 | 26 | 1 | 18 | 11 | |
| 55 | 97 | 69 | 41 | 39 | 32 | 27 | 16 | 11 | 18 | 26 | 1 | 18 | 10 | |
| 56 | 94 | 66 | 41 | 38 | 32 | 26 | 14 | 9 | 18 | 26 | 1 | 18 | 10 | |
| 57 | 94 | 64 | 39 | 38 | 30 | 26 | 11 | 7 | 17 | 25 | 1 | 18 | 9 | |

| | | | | | | | | | | | | | |
|---|---|---|---|---|---|---|---|---|---|---|---|---|---|
| 58 | 93 | 63 | 39 | 37 | 30 | 25 | 9 | 6 | 17 | 25 | | 18 | 8 | |
| 59 | 93 | 62 | 36 | 37 | 29 | 23 | 9 | 6 | 17 | 25 | | 17 | 8 | |
| 60 | 92 | 62 | 35 | 36 | 29 | 23 | 9 | 6 | 16 | 25 | | 17 | 8 | |
| 61 | 91 | **62** | 34 | 36 | 28 | 23 | 8 | 5 | 16 | 25 | | 17 | 8 | |

**Table A1** - continued

| | | | | | | Authors | | | | | | | | |
|---|---|---|---|---|---|---|---|---|---|---|---|---|---|---|
| | L. Leyde-sdorff | W. Glänzel | H. Moed | A. van Raan | R. Rousseau | A. Schubert | B. Martin | F. Narin | E. Garfield | T. Braun | H. Small | L. Egghe | P. Ingwe-Rsen | H.D. White |
| 62 | 90 | 59 | 33 | 35 | 28 | 23 | 7 | 4 | 15 | 24 | | 17 | 7 | |
| 63 | 89 | 59 | 32 | 34 | 28 | 22 | 7 | 2 | 15 | 24 | | 16 | 7 | |
| 64 | 87 | 57 | 30 | 32 | 27 | 22 | 7 | 1 | 14 | 24 | | 16 | 6 | |
| 65 | 87 | 57 | 29 | 31 | 27 | 21 | 7 | 0 | 14 | 23 | | 16 | 5 | |
| 66 | 87 | 56 | 29 | 31 | 27 | 20 | 6 | 0 | 14 | 23 | | 16 | 5 | |
| 67 | 87 | 55 | 28 | 30 | 27 | 19 | 6 | 0 | 14 | 23 | | 16 | 4 | |
| 68 | 86 | 54 | 27 | 29 | 27 | **17** | 5 | 0 | 14 | 22 | | 15 | 4 | |
| 69 | 83 | 54 | 27 | 29 | 26 | 17 | 5 | | 13 | 22 | | 14 | 4 | |
| 70 | 83 | 53 | 27 | 28 | 26 | 16 | 5 | | 13 | 21 | | 14 | 4 | |
| 71 | 83 | 52 | 25 | 28 | 25 | 15 | 5 | | 12 | 21 | | 14 | 4 | |
| 72 | 82 | 52 | 25 | 28 | 24 | 15 | 4 | | 12 | 21 | | 14 | 3 | |
| 73 | 82 | 52 | 23 | 27 | 24 | 15 | 4 | | 12 | 21 | | 14 | 3 | |
| 74 | 82 | 52 | 23 | 25 | 24 | 15 | 4 | | 12 | 20 | | 14 | 3 | |
| 75 | 80 | 51 | 22 | 24 | 24 | 14 | 4 | | 12 | 19 | | 14 | 3 | |
| 76 | 80 | 49 | 22 | 24 | 23 | 14 | 4 | | 12 | 19 | | 13 | 3 | |
| 77 | 80 | 48 | 22 | 24 | 23 | 14 | 3 | | 10 | 19 | | 13 | 2 | |
| 78 | 79 | 48 | 19 | 23 | 23 | 14 | 3 | | 10 | 18 | | 13 | 2 | |
| 79 | 79 | 48 | 18 | 22 | 23 | 13 | 2 | | 10 | 18 | | 13 | 2 | |
| 80 | 78 | 47 | 16 | 21 | 23 | 13 | 2 | | 9 | 17 | | 12 | 2 | |
| 81 | 78 | 46 | 16 | 21 | 23 | 13 | 1 | | 9 | 17 | | 12 | 1 | |
| 82 | 77 | 45 | 15 | 20 | 23 | 13 | 1 | | 9 | 16 | | 12 | 1 | |
| 83 | 77 | 45 | 14 | 20 | 22 | 13 | 1 | | 9 | 16 | | 11 | 1 | |
| 84 | 77 | 44 | 13 | 20 | 22 | 12 | 1 | | 8 | 16 | | 11 | 1 | |
| 85 | 76 | 44 | 12 | 20 | 21 | 12 | 1 | | 7 | 16 | | 11 | 1 | |
| 86 | 76 | 43 | 12 | 18 | 21 | 12 | | | **7** | 16 | | 11 | 1 | |
| 87 | 76 | 43 | 11 | 18 | 21 | 12 | | | 7 | 16 | | 11 | 1 | |
| 88 | 75 | 42 | 11 | 17 | 21 | 10 | | | 6 | 16 | | 11 | 1 | |
| 89 | 75 | 41 | 10 | 17 | 20 | 10 | | | 6 | 15 | | 11 | | |
| 90 | 74 | 40 | 10 | 17 | 20 | 9 | | | 6 | 15 | | 10 | | |
| 91 | 74 | 38 | 9 | 17 | 20 | 9 | | | 6 | 15 | | 10 | | |
| 92 | 73 | 38 | 9 | 13 | 20 | 9 | | | 6 | 15 | | 10 | | |
| 93 | 71 | 38 | 9 | 13 | 20 | 9 | | | 5 | 15 | | 10 | | |
| 94 | 70 | 36 | 8 | 13 | 19 | 8 | | | 5 | 15 | | 10 | | |
| 95 | 69 | 36 | 8 | 13 | 19 | 8 | | | 5 | 15 | | 10 | | |
| 96 | 68 | 36 | 8 | 12 | 18 | 8 | | | 5 | 14 | | 10 | | |
| 97 | 68 | 35 | 7 | 11 | 18 | 7 | | | 5 | 14 | | 10 | | |
| 98 | 68 | 34 | 7 | 10 | 18 | 7 | | | 5 | 14 | | 10 | | |
| 99 | 67 | 34 | 7 | 10 | 18 | 6 | | | 5 | 14 | | 10 | | |
| 100 | 66 | 34 | 7 | 10 | 17 | 6 | | | 5 | 14 | | 10 | | |
| 101 | 65 | 33 | 7 | 9 | 17 | 6 | | | 5 | 14 | | 9 | | |
| 102 | 64 | 33 | 6 | 9 | 17 | 6 | | | 5 | 14 | | 9 | | |
| 103 | 64 | 33 | 6 | 9 | 17 | 6 | | | 5 | 14 | | 9 | | |
| 104 | 63 | 32 | 6 | 9 | 17 | 6 | | | 4 | 14 | | 9 | | |
| 105 | 63 | 32 | 5 | 8 | 17 | 5 | | | 4 | 13 | | 9 | | |
| 106 | 62 | 32 | 5 | 8 | 16 | 5 | | | 4 | 13 | | 9 | | |
| 107 | 61 | 31 | 5 | 8 | 16 | 5 | | | | 13 | | 9 | | |
| 108 | 61 | 31 | 4 | 8 | 16 | 4 | | | | 13 | | 9 | | |
| 109 | 60 | 31 | 4 | 7 | 16 | 4 | | | | 13 | | 9 | | |
| 110 | 60 | 30 | 4 | 7 | 16 | 4 | | | | 13 | | 9 | | |
| 111 | 59 | 30 | 4 | 6 | 16 | 4 | | | | 13 | | 9 | | |
| 112 | 59 | 30 | 4 | 5 | 15 | 4 | | | | 13 | | 9 | | |
| 113 | 58 | 29 | 3 | 5 | 15 | 3 | | | | 13 | | 8 | | |
| 114 | 58 | 28 | 3 | 4 | 15 | 3 | | | | 12 | | 8 | | |
| 115 | 56 | 27 | 3 | 4 | 15 | 3 | | | | 12 | | 8 | | |
| 116 | 56 | 27 | 2 | 4 | 14 | 3 | | | | 12 | | 8 | | |
| 117 | 56 | 27 | 2 | 3 | 14 | 2 | | | | 12 | | 8 | | |

| | 118 | 55 | 27 | 2 | 3 | 14 | 2 | | | | 12 | | 8 | | |
|---|---|---|---|---|---|---|---|---|---|---|---|---|---|---|---|
| | 119 | 55 | 27 | 1 | 3 | 14 | 2 | | | | 12 | | 7 | | |
| | 120 | 55 | 27 | 1 | 2 | 14 | 2 | | | | 11 | | 7 | | |
| | 121 | 55 | 27 | 1 | 2 | 13 | 2 | | | | 11 | | 7 | | |
| | 122 | 53 | 26 | 1 | 1 | 13 | 2 | | | | 10 | | 7 | | |

**Table A1** - continued

| | Authors | | | | | | | | | | | | | |
|---|---|---|---|---|---|---|---|---|---|---|---|---|---|---|
| | L. Leydesdorff | W. Glänzel | H. Moed | A. Van Raan | R. Rousseau | A. Schubert | B. Martin | F. Narin | E. Garfield | T. Braun | H. Small | L. Egghe | P. Ingwe-Rsen | H.D. White |
| 123 | 52 | 26 | 1 | 1 | 13 | 2 | | | | 10 | | 7 | | |
| 124 | 51 | 26 | 1 | 1 | 13 | 2 | | | | 10 | | 7 | | |
| 125 | 51 | 26 | 1 | | 13 | 2 | | | | 10 | | 6 | | |
| 126 | 49 | 26 | 1 | | 13 | 2 | | | | 10 | | 6 | | |
| 127 | 49 | 25 | 1 | | 13 | 2 | | | | 9 | | 6 | | |
| 128 | 48 | 25 | | | 13 | 2 | | | | 9 | | 6 | | |
| 129 | 47 | 24 | | | 13 | 1 | | | | 9 | | 6 | | |
| 130 | 46 | 24 | | | 13 | 1 | | | | 9 | | 6 | | |
| 131 | 44 | 24 | | | 13 | 1 | | | | 9 | | 5 | | |
| 132 | 44 | 23 | | | 13 | 1 | | | | 9 | | 5 | | |
| 133 | 44 | 23 | | | 12 | 1 | | | | 9 | | 5 | | |
| 134 | 44 | 23 | | | 12 | 1 | | | | 9 | | 5 | | |
| 135 | 43 | 23 | | | 12 | 1 | | | | 9 | | 5 | | |
| 136 | 43 | 23 | | | 12 | 1 | | | | 8 | | 5 | | |
| 137 | 43 | 22 | | | 11 | 1 | | | | 8 | | 5 | | |
| 138 | 43 | 21 | | | 11 | 1 | | | | 8 | | 5 | | |
| 139 | 43 | 21 | | | 11 | 1 | | | | 8 | | 5 | | |
| 140 | 43 | 21 | | | 11 | 1 | | | | 8 | | 5 | | |
| 141 | 42 | 20 | | | 11 | 1 | | | | 8 | | 5 | | |
| 142 | 42 | 19 | | | 11 | 1 | | | | 8 | | 5 | | |
| 143 | 41 | 19 | | | 11 | | | | | 8 | | 5 | | |
| 144 | 41 | 18 | | | 11 | | | | | 8 | | 4 | | |
| 145 | 41 | 18 | | | 10 | | | | | 7 | | 4 | | |

Note (citations with the corresponding indices are colored):

| Bold | the $h$-index |
|---|---|
| Blue | the $g$-index |
| Orange | the $har$-index |

We present below the calculations via software *Wolfram Mathematica* the values $S_h := \sum_{j=1}^{h}\frac{1}{cit_j}$,

$S_g := \sum_{j=1}^{g}\frac{1}{cit_j}$, $S_{har+1} := \sum_{j=1}^{har+1}\frac{1}{cit_j}$ and $S_{har} := \sum_{j=1}^{har}\frac{1}{cit_j}$ from Table 1 of all 14 Price awardees, respectively.

**L. Leydesdorf**
1/
1/3964+1/647+1/544+1/479+1/424+1/422+1/410+1/401+1/376+1/352+1/352+1/301+1/291+1/266+1/252+1/242+1/241+1/232+1/225+1/221+1/184+1/179+1/171+1/166+1/165+1/157+1/154+1/147+1/147+1/147+1/141+1/138+1/135+1/133+1/133+1/129+1/129+1/127+1/127+1/126+1/125+1/125+1/125+1/125+1/124+1/123+1/122+1/122+1/118+1/108+1/106+1/102+1/100+1/98+1/97+1/94+1/94+1/93+1/93+1/92+1/91+1/90+1/89+1/87+1/87+1/87+1/87+1/86+1/83+1/83+1/83+1/82+1/82+1/82+1/80+1/80+1/80+1/79+1/79+1/78+1/78+1/77+1/77+1/77+1/76+1/76+1/76+1/75+1/75+1/74+1/74+1/73+1/71+1/70+1/69+1/68+1/68+1/68+1/67+1/66+1/65+1/64+1/64+1/63+1/63+1/62+1/61+1/61+1/60+1/60+1/59+1/59+1/58+1/58+1/56+1/56+1/56+1/55+1/55+1/55+1/55+1/53+1/52+1/51+1/51+1/49+1/49+1/48+1/47+1/46+1/44+1/44+1/44+1/44+1/43+1/43+1/43+1/43+1/43+1/43+1/42+1/42+1/41+1/41+1/41.0
1.78723

1/3964+1/647+1/544+1/479+1/424+1/422+1/410+1/401+1/376+1/352+1/352+1/301+1/291+1/266+1/252+1/242+1/241+1/232+1/225+1/221+1/184+1/179+1/171+1/166+1/165+1/157+1/154+1/147+1/147+1/147+1/141+1/138+1/135+1/133+1/133+1/129+1/129+1/127+1/127+1/126+1/125+1/125+1/125+1/125+1/124+1/123+1/122+1/122+1/118+1/108+1/106+1/102+1/100+1/98+1/97+1/94+1/94+1/93+1/93+1/92+1/91+1/90+1/89+1/87+1/87+1/87+1/87+1/86+1/83+1/83+1/83+1/82+1/82+1/82+1/80+1/80+1/80+1/79+1/79.0
0.610183

1/3964+1/647+1/544+1/479+1/424+1/422+1/410+1/401+1/376+1/352+1/352+1/301+1/291+1/266+1/252+1/242+1/241+1/232+1/225+1/221+1/184+1/179+1/171+1/166+1/165+1/157+1/154+1/147+1/147+1/147+1/141+1/138+1/135+1/133+1/133+1/129+1/129+1/127+1/127+1/126+1/125+1/125+1/125+1/125+1/124+1/123+1/122+1/122+1/118+1/108+1/106+1/102+1/100+1/98+1/97+1/94+1/94+1/93+1/93+1/92+1/91+1/90+1/89+1/87+1/87+1/87+1/87+1/86+1/83+1/83+1/83+1/82+1/82+1/82+1/80+1/80+1/80+1/79+1/79+1/78+1/78+1/77+1/77+1/77+1/76+1/76+1/76+1/75+1/75+1/74+1/74+1/73+1/71+1/70+1/69+1/68+1/68+1/68+1/67+1/66+1/65+1/64+1/64+1/63+1/63+1/62.0+1/61.0
1.00961

1/3964+1/647+1/544+1/479+1/424+1/422+1/410+1/401+1/376+1/352+1/352+1/301+1/291+1/266+1/252+1/242+1/241+1/232+1/225+1/221+1/184+1/179+1/171+1/166+1/165+1/157+1/154+1/147+1/147+1/147+1/141+1/138+1/135+1/133+1/133+1/129+1/129+1/127+1/127+1/126+1/125+1/125+1/125+1/125+1/124+1/123+1/122+1/122+1/118+1/108+1/106+1/102+1/100+1/98+1/97+1/94+1/94+1/93+1/93+1/92+1/91+1/90+1/89+1/87+1/87+1/87+1/87+1/86+1/83+1/83+1/83+1/82+1/82+1/82+1/80+1/80+1/80+1/79+1/79+1/78+1/78+1/77+1/77+1/77+1/76+1/76+1/76+1/75+1/75+1/74+1/74+1/73+1/71+1/70+1/69+1/68+1/68+1/68+1/67+1/66+1/65+1/64+1/64+1/63+1/63+1/62.0
0.993218

**W. Glänzel**

1/534+1/449+1/417+1/364+1/273+1/255+1/249+1/220+1/211+1/201+1/188+1/163+1/159+1/155+1/147+1/145+1/141+1/135+1/133+1/127+1/126+1/126+1/121+1/111+1/109+1/108+1/96+1/94+1/92+1/89+1/89+1/89+1/85+1/85+1/84+1/84+1/83+1/80+1/80+1/78+1/76+1/76+1/75+1/75+1/75+1/74+1/73+1/73+1/73+1/72+1/72+1/71+1/71+1/70+1/69+1/66+1/64+1/63+1/62+1/62+1/62.0
0.618

1/534+1/449+1/417+1/364+1/273+1/255+1/249+1/220+1/211+1/201+1/188+1/163+1/159+1/155+1/147+1/145+1/141+1/135+1/133+1/127+1/126+1/126+1/121+1/111+1/109+1/108+1/96+1/94+1/92+1/89+1/89+1/89+1/85+1/85+1/84+1/84+1/83+1/80+1/80+1/78+1/76+1/76+1/75+1/75+1/75+1/74+1/73+1/73+1/73+1/72+1/72+1/71+1/71+1/70+1/69+1/66+1/64+1/63+1/62+1/62+1/62+1/59+1/59+1/57+1/57+1/56+1/55+1/54+1/54+1/53+1/52+1/52+1/52+1/52+1/51+1/49+1/48+1/48+1/48+1/47+1/46+1/45+1/45+1/44+1/44+1/43+1/43+1/42+1/41+1/40+1/38+1/38+1/38+1/36+1/36+1/36+1/35+1/34+1/34.0
1.46067

1/534+1/449+1/417+1/364+1/273+1/255+1/249+1/220+1/211+1/201+1/188+1/163+1/159+1/155+1/147+1/145+1/141+1/135+1/133+1/127+1/126+1/126+1/121+1/111+1/109+1/108+1/96+1/94+1/92+1/89+1/89+1/89+1/85+1/85+1/84+1/84+1/83+1/80+1/80+1/78+1/76+1/76+1/75+1/75+1/75+1/74+1/73+1/73+1/73+1/72+1/72+1/71+1/71+1/70+1/69+1/66+1/64+1/63+1/62+1/62+1/62+1/59+1/59+1/57+1/57+1/56+1/55+1/54+1/54+1/53+1/52+1/52+1/52+1/52+1/51+1/49+1/48+1/48+1/48+1/47+1/46.0
1.00138

1/534+1/449+1/417+1/364+1/273+1/255+1/249+1/220+1/211+1/201+1/188+1/163+1/159+1/155+1/147+1/145+1/141+1/135+1/133+1/127+1/126+1/126+1/121+1/111+1/109+1/108+1/96+1/94+1/92+1/89+1/89+1/89+1/85+1/85+1/84+1/84+1/83+1/80+1/80+1/78+1/76+1/76+1/75+1/75+1/75+1/74+1/73+1/73+1/73+1/72+1/72+1/71+1/71+1/70+1/69+1/66+1/64+1/63+1/62+1/62+1/62+1/59+1/59+1/57+1/57+1/56+1/55+1/54+1/54+1/53+1/52+1/52+1/52+1/52+1/51+1/49+1/48+1/48+1/48+1/47.0
0.979645

**H.F. Moed**

1/417+1/417+1/380+1/350+1/303+1/290+1/288+1/248+1/197+1/183+1/182+1/174+1/166+1/152+1/146+1/122+1/118+1/118+1/117+1/103+1/95+1/92+1/88+1/80+1/76+1/75+1/75+1/73+1/71+1/70+1/69+1/69+1/63+1/62+1/61+1/60+1/60+1/58+1/57+1/57+1/56+1/55+1/53+1/53+1/52+1/52+1/50+1/49+1/49.0
0.572423

1/417+1/417+1/380+1/350+1/303+1/290+1/288+1/248+1/197+1/183+1/182+1/174+1/166+1/152+1/146+1/122+1/118+1/118+1/117+1/103+1/95+1/92+1/88+1/80+1/76+1/75+1/75+1/73+1/71+1/70+1/69+1/69+1/63+1/62+1/61+1/60+1/60+1/58+1/57+1/57+1/56+1/55+1/53+1/53+1/52+1/52+1/50+1/49+1/49.0+1/47+1/45+1/44+1/43+1/43+1/41+1/41+1/39+1/39+1/36+1/35+1/34+1/33+1/32+1/30+1/29+1/29+1/28+1/27+1/27+1/27+1/25+1/25+1/23+1/23+1/22+1/22+1/22+1/19+1/18+1/16+1/16+1/15+1/14+1/13+1/12+1/11.0
2.10743

1/417+1/417+1/380+1/350+1/303+1/290+1/288+1/248+1/197+1/183+1/182+1/174+1/166+1/152+1/146+1/122+1/118+1/118+1/117+1/103+1/95+1/92+1/88+1/80+1/76+1/75+1/75+1/73+1/71+1/70+1/69+1/69+1/63+1/62+1/61+1/60+1/60+1/58+1/57+1/57+1/56+1/55+1/53+1/53+1/52+1/52+1/50+1/49+1/49.0+1/47+1/45+1/44+1/43+1/43+1/41+1/41+1/39+1/39+1/36+1/35+1/34+1/33+1/32+1/30+1/29
1.00035

1/417+1/417+1/380+1/350+1/303+1/290+1/288+1/248+1/197+1/183+1/182+1/174+1/166+1/152+1/146+1/122+1/118+1/118+1/117+1/103+1/95+1/92+1/88+1/80+1/76+1/75+1/75+1/73+1/71+1/70+1/69+1/69+1/63+1/62+1/61+1/60+1/60+1/58+1/57+1/57+1/56+1/55+1/53+1/53+1/52+1/52+1/50+1/49+1/49.0+1/47+1/45+1/44+1/43+1/43+1/41+1/41+1/39+1/39+1/36+1/35+1/34+1/33+1/32+1/30
0.965871

**A. Van Raan**

1/522+1/455+1/381+1/350+1/303+1/290+1/284+1/283+1/276+1/206+1/172+1/166+1/157+1/152+1/150+1/146+1/142+1/136+1/130+1/111+1/110+1/103+1/98+1/94+1/92+1/86+1/82+1/81+1/81+1/77+1/76+1/75+1/73+1/72+1/72+1/68+1/65+1/63+1/62+1/60+1/59+1/59+1/59+1/55+1/52+1/51+1/48+1/48.0
0.509549

1/522+1/455+1/381+1/350+1/303+1/290+1/284+1/283+1/276+1/206+1/172+1/166+1/157+1/152+1/150+1/146+1/142+1/136+1/130+1/111+1/110+1/103+1/98+1/94+1/92+1/86+1/82+1/81+1/81+1/77+1/76+1/75+1/73+1/72+1/72+1/68+1/65+1/63+1/62+1/60+1/59+1/59+1/59+1/55+1/52+1/51+1/48+1/48.0+1/47+1/45+1/41+1/41+1/39+1/39+1/39+1/38+1/38+1/37+1/37+1/36+1/36+1/35+1/35+1/34+1/32+1/31+1/31+1/30+1/29+1/29+1/28+1/28+1/28+1/27+1/25+1/24+1/24+1/24+1/23+1/22+1/21+1/21+1/20+1/20+1/20+1/20+1/18+1/18+1/17+1/17.0
2.04772

1/522+1/455+1/381+1/350+1/303+1/290+1/284+1/283+1/276+1/206+1/172+1/166+1/157+1/152+1/150+1/146+1/142+1/136+1/130+1/111+1/110+1/103+1/98+1/94+1/92+1/86+1/82+1/81+1/81+1/77+1/76+1/75+1/73+1/72+1/72+1/68+1/65+1/63+1/62+1/60+1/59+1/59+1/59+1/55+1/52+1/51+1/48+1/48.0+1/47+1/45+1/41+1/41+1/39+1/39+1/39+1/38+1/38+1/37+1/37+1/36+1/36+1/35+1/35+1/34+1/32+1/31+1/31
1.02331

1/522+1/455+1/381+1/350+1/303+1/290+1/284+1/283+1/276+1/206+1/172+1/166+1/157+1/152+1/150+1/146+1/142+1/136+1/130+1/111+1/110+1/103+1/98+1/94+1/92+1/86+1/82+1/81+1/81+1/77+1/76+1/75+1/73+1/72+1/72+1/68+1/65+1/63+1/62+1/60+1/59+1/59+1/59+1/55+1/52+1/51+1/48+1/48.0+1/47+1/45+1/41+1/41+1/39+1/39+1/39+1/38+1/38+1/37+1/37+1/36+1/36+1/35+1/35+1/34+1/32+1/31
0.991056

**R. Rousseau**

1/1050+1/484+1/478+1/265+1/242+1/175+1/153+1/126+1/118+1/111+1/110+1/97+1/96+1/88+1/87+1/75+1/74+1/72+1/70+1/67+1/66+1/64+1/63+1/61+1/58+1/56+1/56+1/55+1/51+1/49+1/48+1/48+1/47+1/46+1/45+1/45+1/45+1/45+1/44+1/44+1/43+1/43+1/43.0
0.636298

1/1050+1/484+1/478+1/265+1/242+1/175+1/153+1/126+1/118+1/111+1/110+1/97+1/96+1/88+1/87+1/75+1/74+1/72+1/70+1/67+1/66+1/64+1/63+1/61+1/58+1/56+1/56+1/55+1/51+1/49+1/48+1/48+1/47+1/46+1/45+1/45+1/45+1/45+1/44+1/44+1/43+1/43+1/43.0+1/42+1/42+1/42+1/41+1/40+1/39+1/39+1/39+1/37+1/37+1/33+1/32+1/32+1/30+1/30+1/29+1/29+1/28+1/28+1/28+1/27+1/27+1/27+1/27+1/27+1/26+1/26+1/25+1/24+1/24+1/24+1/24+1/23+1/23+1/23+1/23.0
1.86636

1/1050+1/484+1/478+1/265+1/242+1/175+1/153+1/126+1/118+1/111+1/110+1/97+1/96+1/88+1/87+1/75+1/74+1/72+1/70+1/67+1/66+1/64+1/63+1/61+1/58+1/56+1/56+1/55+1/51+1/49+1/48+1/48+1/47+1/46+1/45+1/45+1/45+1/45+1/44+1/44+1/43+1/43+1/43.0+1/42+1/42+1/42+1/41+1/40+1/39+1/39+1/39+1/37+1/37+1/33+1/32+1/32+1/30
1.01423

1/1050+1/484+1/478+1/265+1/242+1/175+1/153+1/126+1/118+1/111+1/110+1/97+1/96+1/88+1/87+1/75+1/74+1/72+1/70+1/67+1/66+1/64+1/63+1/61+1/58+1/56+1/56+1/55+1/51+1/49+1/48+1/48+1/47+1/46+1/45+1/45+1/45+1/45+1/44+1/44+1/43+1/43+1/43.0+1/42+1/42+1/42+1/41+1/40+1/39+1/39+1/39+1/37+1/37+1/33+1/32
0.980897

**A. Schubert**

1/1858+1/449+1/362+1/273+1/255+1/211+1/159+1/158+1/155+1/151+1/147+1/133+1/127+1/120+1/119+1/105+1/98+1/92+1/89+1/75+1/73+1/73+1/72+1/71+1/70+1/69+1/69+1/64+1/63+1/62+1/62+1/57+1/53+1/52+1/47+1/46+1/46+1/45+1/45+1/43+1/42+1/42.0
0.535758

(1/1858+1/449+1/362+1/273+1/255+1/211+1/159+1/158+1/155+1/151+1/147+1/133+1/127+1/120+1/119+1/105+1/98+1/92+1/89+1/75+1/73+1/73+1/72+1/71+1/70+1/69+1/69+1/64+1/63+1/62+1/62+1/57+1/53+1/52+1/47+1/46+1/46+1/45+1/45+1/43+1/42+1/42.0)+1/40+1/37+1/36+1/36+1/36+1/33+1/32+1/31+1/31+1/27+1/27+1/27+1/27+1/26+1/26+1/25+1/23+1/23+1/23+1/23+1/22+1/22+1/21+1/20+1/19+1/17+1/17+1/16+1/15+1/15+1/15+1/15+1/14+1/14+1/14+1/14+1/13+1/13+1/13+1/13+1/13+1/12+1/12.0
2.76114

(1/1858+1/449+1/362+1/273+1/255+1/211+1/159+1/158+1/155+1/151+1/147+1/133+1/127+1/120+1/119+1/105+1/98+1/92+1/89+1/75+1/73+1/73+1/72+1/71+1/70+1/69+1/69+1/64+1/63+1/62+1/62+1/57+1/53+1/52+1/47+1/46+1/46+1/45+1/45+1/43+1/42+1/42.0)+1/40+1/37+1/36+1/36+1/36+1/33+1/32+1/31+1/31+1/27+1/27+1/27+1/27+1/26+1/26
1.02226

(1/1858+1/449+1/362+1/273+1/255+1/211+1/159+1/158+1/155+1/151+1/147+1/133+1/127+1/120+1/119+1/105+1/98+1/92+1/89+1/75+1/73+1/73+1/72+1/71+1/70+1/69+1/69+1/64+1/63+1/62+1/62+1/57+1/53+1/52+1/47+1/46+1/46+1/45+1/45+1/43+1/42+1/42.0)+1/40+1/37+1/36+1/36+1/36+1/33+1/32+1/31+1/31+1/27+1/27+1/27+1/27+1/26
0.983797

1/40+1/37+1/36+1/36+1/36+1/33+1/32+1/31+1/31+1/27+1/27+1/27+1/27+1/26+1/26+1/25+1/23+1/23+1/23+1/23+1/22+1/22+1/21+1/20+1/19+1/17+1/17+1/16+1/15+1/15+1/15+1/15+1/14+1/14+1/14+1/14+1/13+1/13+1/13+1/13+1/13+1/12+1/12.0
2.22538

**B. Martin**

1/1828+1/728+1/411+1/391+1/371+1/362+1/282+1/244+1/223+1/196+1/179+1/131+1/128+1/122+1/114+1/101+1/77+1/73+1/72+1/69+1/68+1/66+1/65+1/64+1/63+1/58+1/57+1/55+1/50+1/48+1/48+1/47+1/46+1/46+1/45+1/40+1/40+1/40.0
0.48592

1/1828+1/728+1/411+1/391+1/371+1/362+1/282+1/244+1/223+1/196+1/179+1/131+1/128+1/122+1/114+1/101+1/77+1/73+1/72+1/69+1/68+1/66+1/65+1/64+1/63+1/58+1/57+1/55+1/50+1/48+1/48+1/47+1/46+1/46+1/45+1/40+1/40+1/40.0+1/32+1/32+1/31+1/28+1/24+1/23+1/23+1/23+1/22+1/21+1/21+1/21+1/20+1/20+1/20+1/17+1/16+1/14+1/11+1/9+1/9+1/9+1/8+1/7+1/7+1/7+1/7+1/6+1/6+1/5+1/5+1/5+1/5+1/4+1/4+1/4+1/4+1/4+1/3+1/3+1/2+1/2+1/1+1/1+1/1+1/1+1/1+1/1
11.4902

1/1828+1/728+1/411+1/391+1/371+1/362+1/282+1/244+1/223+1/196+1/179+1/131+1/128+1/122+1/114+1/101+1/77+1/73+1/72+1/69+1/68+1/66+1/65+1/64+1/63+1/58+1/57+1/55+1/50+1/48+1/48+1/47+1/46+1/46+1/45+1/40+1/40+1/40.0+1/32+1/32+1/31+1/28+1/24+1/23+1/23+1/23+1/22+1/21+1/21+1/21+1/20
1.02681

1/1828+1/728+1/411+1/391+1/371+1/362+1/282+1/244+1/223+1/196+1/179+1/131+1/128+1/122+1/114+1/101+1/77+1/73+1/72+1/69+1/68+1/66+1/65+1/64+1/63+1/58+1/57+1/55+1/50+1/48+1/48+1/47+1/46+1/46+1/45+1/40+1/40+1/40.0+1/32+1/32+1/31+1/28+1/24+1/23+1/23+1/23+1/22+1/21+1/21+1/21
0.976806

**F. Narin**

(1/803+1/780+1/519+1/481+1/443+1/354+1/339+1/306+1/283+1/226+1/193+1/183+1/159+1/151+1/138+1/121+1/114+1/94+1/87+1/81+1/76+1/73+1/73+1/68+1/67+1/61+1/61+1/55+1/54+1/50+1/43+1/43+1/43+1/42+1/41+1/40+1/40+1/38.0)
0.462067

1/803+1/780+1/519+1/481+1/443+1/354+1/339+1/306+1/283+1/226+1/193+1/183+1/159+1/151+1/138+1/121+1/114+1/94+1/87+1/81+1/76+1/73+1/73+1/68+1/67+1/61+1/61+1/55+1/54+1/50+1/43+1/43+1/43+1/42+1/41+1/40+1/40+1/38+1/32+1/32+1/31+1/28+1/24+1/23+1/23+1/23+1/22+1/21+1/21+1/21+1/20+1/20+1/20+1/17+1/16+1/14+1/11+1/9+1/9+1/9+1/8+1/7+1/7+1/7+1/7+1/6+1/6+1/5.0
2.94971

(1/803+1/780+1/519+1/481+1/443+1/354+1/339+1/306+1/283+1/226+1/193+1/183+1/159+1/151+1/138+1/121+1/114+1/94+1/87+1/81+1/76+1/73+1/73+1/68+1/67+1/61+1/61+1/55+1/54+1/50+1/43+1/43+1/43+1/42+1/41+1/40+1/40+1/38.0)+1/37+1/34+1/32+1/27+1/24+1/21+1/21+1/20+1/18+1/18+1/17+1/14
1.01506

(1/803+1/780+1/519+1/481+1/443+1/354+1/339+1/306+1/283+1/226+1/193+1/183+1/159+1/151+1/138+1/121+1/114+1/94+1/87+1/81+1/76+1/73+1/73+1/68+1/67+1/61+1/61+1/55+1/54+1/50+1/43+1/43+1/43+1/42+1/41+1/40+1/40+1/38.0)+1/37+1/34+1/32+1/27+1/24+1/21+1/21+1/20+1/18+1/18+1/17
0.943

**E. Garfield**

1/1888+1/1711+1/1635+1/597+1/589+1/430+1/264+1/255+1/252+1/238+1/238+1/178+1/155+1/154+1/153+1/151+1/139+1/132+1/114+1/113+1/111+1/91+1/88+1/83+1/74+1/73+1/69+1/65+1/64+1/63+1/63+1/53+1/51+1/49+1/45+1/42+1/39.0
0.370045

1/1888+1/1711+1/1635+1/597+1/589+1/430+1/264+1/255+1/252+1/238+1/238+1/178+1/155+1/154+1/153+1/151+1/139+1/132+1/114+1/113+1/111+1/91+1/88+1/83+1/74+1/73+1/69+1/65+1/64+1/63+1/63+1/53+1/51+1/49+1/45+1/42+1/39+1/30+1/28+1/27+1/27+1/25+1/23+1/22+1/22+1/21+1/21+1/21+1/21+1/21+1/20+1/20+1/20+1/19+1/18+1/18+1/17+1/17+1/16+1/15+1/15+1/14+1/14+1/14+1/14+1/13+1/13+1/12+1/12+1/12+1/12+1/12+1/12+1/10+1/10+1/10+1/9+1/9+1/9+1/9+1/8+1/7+1/7+1/7+1/6+1/6+1/6+1/6+1/6+1/5+1/5+1/5+1/5+1/5+1/5+1/5+1/5+1/5+1/5+1/4+1/4+1/4.0
7.57378

1/1888+1/1711+1/1635+1/597+1/589+1/430+1/264+1/255+1/252+1/238+1/238+1/178+1/155+1/154+1/153+1/151+1/139+1/132+1/114+1/113+1/111+1/91+1/88+1/83+1/74+1/73+1/69+1/65+1/64+1/63+1/63+1/53+1/51+1/49+1/45+1/42+1/39+1/30+1/28+1/27+1/27+1/25+1/23+1/22+1/22+1/21+1/21+1/21+1/21+1/21+1/20+1/20.0
1.02565

1/1888+1/1711+1/1635+1/597+1/589+1/430+1/264+1/255+1/252+1/238+1/238+1/178+1/155+1/154+1/153+1/151+1/139+1/132+1/114+1/113+1/111+1/91+1/88+1/83+1/74+1/73+1/69+1/65+1/64+1/63+1/63+1/53+1/51+1/49+1/45+1/42+1/39+1/30+1/28+1/27+1/27+1/25+1/23+1/22+1/22+1/21+1/21+1/21+1/21+1/21+1/20.0
0.97565

**T. Braun**

(1/449+1/362+1/255+1/158+1/151+1/147+1/138+1/135+1/127+1/120+1/119+1/96+1/77+1/75+1/73+1/69+1/67+1/61+1/59+1/57+1/55+1/55+1/52+1/52+1/51+1/47+1/46+1/45+1/45+1/45+1/43+1/42+1/40+1/40+1/39+1/37+1/37.0)
0.579546

1/449+1/362+1/255+1/158+1/151+1/147+1/138+1/135+1/127+1/120+1/119+1/96+1/77+1/75+1/73+1/69+1/67+1/61+1/59+1/57+1/55+1/55+1/52+1/52+1/51+1/47+1/46+1/45+1/45+1/45+1/43+1/42+1/40+1/40+1/39+1/37+1/37.0+1/36+1/34+1/33+1/33+1/33+1/32+1/31+1/31+1/30+1/29+1/28+1/28+1/27+1/27+1/27+1/27+1/26+1/26+1/26+1/25+1/25+1/25+1/25+1/25+1/24+1/24+1/24+1/23
1.59467

(1/449+1/362+1/255+1/158+1/151+1/147+1/138+1/135+1/127+1/120+1/119+1/96+1/77+1/75+1/73+1/69+1/67+1/61+1/59+1/57+1/55+1/55+1/52+1/52+1/51+1/47+1/46+1/45+1/45+1/45+1/43+1/42+1/40+1/40+1/39+1/37+1/37.0)+1/36+1/34+1/33+1/33+1/33+1/32+1/31+1/31+1/30+1/29+1/28+1/28+1/27+1/27
1.03673

(1/449+1/362+1/255+1/158+1/151+1/147+1/138+1/135+1/127+1/120+1/119+1/96+1/77+1/75+1/73+1/69+1/67+1/61+1/59+1/57+1/55+1/55+1/52+1/52+1/51+1/47+1/46+1/45+1/45+1/45+1/43+1/42+1/40+1/40+1/39+1/37+1/37.0)+1/36+1/34+1/33+1/33+1/33+1/32+1/31+1/31+1/30+1/29+1/28+1/28+1/27
0.999693

**H. Small**

1/2846+1/592+1/540+1/410+1/242+1/197+1/193+1/190+1/186+1/157+1/154+1/151+1/144+1/129+1/113+1/112+1/111+1/102+1/98+1/90+1/81+1/79+1/71+1/57+1/51+1/49+1/47+1/44+1/43+1/41+1/40+1/39+1/38+1/34.0
0.418116

1/2846+1/592+1/540+1/410+1/242+1/197+1/193+1/190+1/186+1/157+1/154+1/151+1/144+1/129+1/113+1/112+1/111+1/102+1/98+1/90+1/81+1/79+1/71+1/57+1/51+1/49+1/47+1/44+1/43+1/41+1/40+1/39+1/38+1/34.0+1/21+1/21+1/20+1/19+1/18+1/17+1/17+1/15+1/14+1/10+1/9+1/9+1/6+1/4+1/4+1/4+1/4+1/2+1+1+1+1+1.0
7.91617

1/2846+1/592+1/540+1/410+1/242+1/197+1/193+1/190+1/186+1/157+1/154+1/151+1/144+1/129+1/113+1/112+1/111+1/102+1/98+1/90+1/81+1/79+1/71+1/57+1/51+1/49+1/47+1/44+1/43+1/41+1/40+1/39+1/38+1/34.0+1/21+1/21+1/20+1/19+1/18+1/17+1/17+1/15+1/14+1/10
1.02728

(1/2846+1/592+1/540+1/410+1/242+1/197+1/193+1/190+1/186+1/157+1/154+1/151+1/144+1/129+1/113+1/112+1/111+1/102+1/98+1/90+1/81+1/79+1/71+1/57+1/51+1/49+1/47+1/44+1/43+1/41+1/40+1/39+1/38+1/34.0)+1/21+1/21+1/20+1/19+1/18+1/17+1/17+1/15+1/14
0.927283

**I. Egghe**

1/1455+1/478+1/242+1/206+1/153+1/124+1/110+1/107+1/88+1/81+1/71+1/70+1/70+1/68+1/65+1/50+1/48+1/47+1/46+1/45+1/42+1/41+1/40+1/40+1/39+1/35+1/34+1/34+1/34+1/32.0
0.519228

(1/1455+1/478+1/242+1/206+1/153+1/124+1/110+1/107+1/88+1/81+1/71+1/70+1/70+1/68+1/65+1/50+1/48+1/47+1/46+1/45+1/42+1/41+1/40+1/40+1/39+1/35+1/34+1/34+1/34+1/32.0)+1/30+1/29+1/29+1/28+1/27+1/26+1/25+1/24+1/24+1/24+1/24+1/24+1/23+1/22+1/22+1/22+1/22+1/22+1/20+1/20+1/20+1/19+1/19+1/18+1/18+1/18+1/18+1/18+1/17+1/17+1/17+1/17+1/16+1/16+1/16+1/16+1/16+1/15+1/14.0
2.47075

(1/1455+1/478+1/242+1/206+1/153+1/124+1/110+1/107+1/88+1/81+1/71+1/70+1/70+1/68+1/65+1/50+1/48+1/47+1/46+1/45+1/42+1/41+1/40+1/40+1/39+1/35+1/34+1/34+1/34+1/32.0)+1/30+1/29+1/29+1/28+1/27+1/26+1/25+1/24+1/24+1/24+1/24+1/24+1/23
1.02455

(1/1455+1/478+1/242+1/206+1/153+1/124+1/110+1/107+1/88+1/81+1/71+1/70+1/70+1/68+1/65+1/50+1/48+1/47+1/46+1/45+1/42+1/41+1/40+1/40+1/39+1/35+1/34+1/34+1/34+1/32.0)+1/30+1/29+1/29+1/28+1/27+1/26+1/25+1/24+1/24+1/24+1/24+1/24
0.981073

**P. Iingwersen**

1/482+1/414+1/348+1/232+1/166+1/156+1/155+1/133+1/113+1/70+1/57+1/56+1/49+1/46+1/45+1/42+1/42+1/41+1/39+1/39+1/38+1/35+1/35+1/33+1/29+1/29+1/28.0
0.502717

(1/482+1/414+1/348+1/232+1/166+1/156+1/155+1/133+1/113+1/70+1/57+1/56+1/49+1/46+1/45+1/42+1/42+1/41+1/39+1/39+1/38+1/35+1/35+1/33+1/29+1/29+1/28.0)+1/27+1/27+1/26+1/25+1/24+1/24+1/24+1/23+1/23+1/23+1/22+1/22+1/22+1/18+1/17+1/16+1/16+1/15+1/15+1/14+1/14+1/13+1/13+1/13+1/12+1/12+1/11+1/10+1/10+1/9+1/8+1/8+1/8
2.73708

(1/482+1/414+1/348+1/232+1/166+1/156+1/155+1/133+1/113+1/70+1/57+1/56+1/49+1/46+1/45+1/42+1/42+1/41+1/39+1/39+1/38+1/35+1/35+1/33+1/29+1/29+1/28.0)+1/27+1/27+1/26+1/25+1/24+1/24+1/24+1/23+1/23+1/23+1/22+1/22
1.0016

(1/482+1/414+1/348+1/232+1/166+1/156+1/155+1/133+1/113+1/70+1/57+1/56+1/49+1/46+1/45+1/42+1/42+1/41+1/39+1/39+1/38+1/35+1/35+1/33+1/29+1/29+1/28.0)+1/27+1/27+1/26+1/25+1/24+1/24+1/24+1/23+1/23+1/23+1/22
0.956142

**H.D. White**

1/1029+1/213+1/178+1/148+1/124+1/123+1/104+1/72+1/54+1/50+1/35+1/32+1/29+1/28+1/25+1/24+1/23+1/21+1/20.0
0.449041

1/1029+1/213+1/178+1/148+1/124+1/123+1/104+1/72+1/54+1/50+1/35+1/32+1/29+1/28+1/25+1/24+1/23+1/21+1/20+1/19+1/17+1/15+1/5+1/5+1/2+1/2+1+1.0
4.02716

1/19+1/17+1/15+1/5+1/5.0
1/1029+1/213+1/178+1/148+1/124+1/123+1/104+1/72+1/54+1/50+1/35+1/32+1/29+1/28+1/25+1/24+1/23+1/21+1/20.0+1/19+1/17+1/15+1/5+1/5.0
1.02716

1/1029+1/213+1/178+1/148+1/124+1/123+1/104+1/72+1/54+1/50+1/35+1/32+1/29+1/28+1/25+1/24+1/23+1/21+1/20.0+1/19+1/17+1/15+1/5
0.827163